\documentclass[graybox]{SNmult}

\usepackage{type1cm}        
\usepackage{makeidx}         
\usepackage{graphicx}        
\usepackage{multicol}        
\usepackage[bottom]{footmisc}

\usepackage{newtxtext}       %
\usepackage[varvw]{newtxmath}       

\usepackage{comment}
\newcommand{\vertiii}[1]{{\left\vert\kern-0.25ex\left\vert\kern-0.25ex\left\vert #1 
		\right\vert\kern-0.25ex\right\vert\kern-0.25ex\right\vert}}
\newtheorem{assumption}[theorem]{Assumption}
\newtheorem{problem1}[theorem]{Problem}
\newtheorem{theorem1}[theorem]{Theorem}

\usepackage{pifont}

\makeindex             

\begin{document}
\title*{Reduced Basis Method for Simulating Thermal Transients in Electric Machines}
\author{Herbert Egger\orcidID{0000-0003-3769-8791} 
and\\ Eva-Maria Haslhofer\orcidID{0009-0003-0969-5427}
}
\institute{Herbert Egger \at Johann Radon Institute for Computational and Applied Mathematics, Linz Austria,\\ \email{herbert.egger@ricam.oeaw.ac.at}
\and Eva-Maria Haslhofer \at Institute of Numerical Mathematics, Johannes Kepler University Linz, Austria, \\ \email{eva-maria.haslhofer@jku.at}}

\maketitle
\abstract*{
Reduced Basis (RB) methods provide low-dimensional approximations of parametrized partial differential equations with controllable accuracy. We study the systematic construction of RB methods for transient thermal simulation of rotating electric machines and their calibration to measurements for an induction motor using parameter estimation. The performance of the proposed approach is compared to Finite Element (FE) models and Lumped Parameter Thermal Networks (LPTNs). 
\keywords{Thermal Modeling $\cdot$ Electric Machines $\cdot$ Model Order Reduction}}

\abstract{
Reduced Basis (RB) methods provide low-dimensional approximations of parametrized partial differential equations with controllable accuracy. We discuss the construction of RB approximations for transient thermal simulation of an induction motor and compare their performance to Finite Element (FE) models and Lumped Parameter Thermal Networks (LPTNs) after calibration to measurements. Numerical results demonstrate that RB models can achieve FE accuracy, allow for a fully automatic construction, and offer computational efficiency comparable to traditional LPTNs.
\keywords{Thermal Modeling $\cdot$ Electric Machines $\cdot$ Model Order Reduction}}

\section{Introduction}
\label{haslhofer:sec:1}

Thermal behavior is a key limiting factor in the design and operation of modern electric machines, particularly under increasing power density and efficiency requirements~\cite{haslhofer:US_drive}. Reliable prediction of temperature distributions is essential to ensure safe operation and prevent thermal damage~\cite{haslhofer:laughton}. Since only a limited number of temperature sensors are available in practice, 
state estimation is required to assess the temperature distribution in the machine, and reduced-order models are used for online control applications\cite{haslhofer:hartmann}.

\begin{figure}[ht!]
        \centering
        \vspace*{-2mm}
        \includegraphics[width=0.45\textwidth]{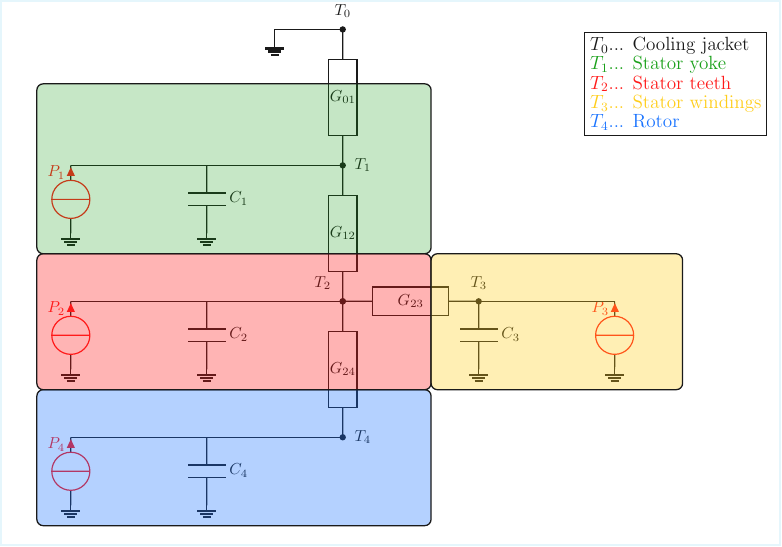}
        \hspace{5mm}
        \includegraphics[ width=0.28\textwidth]{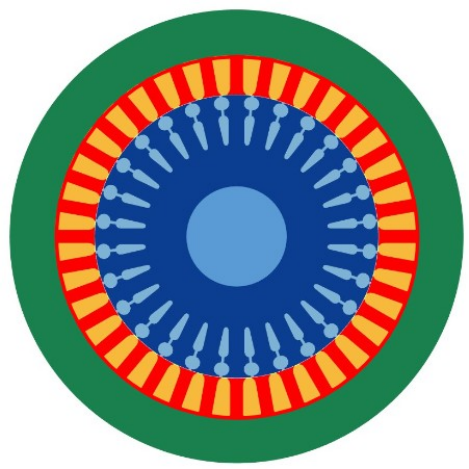}
          \caption{LPTN set-up with five nodes (Left) set up by~\cite{haslhofer:eickhoff}; Cross section of induction machine (Right); Colors indicate affiliation to respective compartments.}%
         \label{haslhofer:fig:LPTN}
\end{figure}

A widely used industrial approach for simulating transients in electric machines are Lumped Parameter Thermal Networks~\cite{haslhofer:staton}. They partition the machine into a small number of thermally homogeneous regions connected by equivalent thermal resistances; see Figure~\ref{haslhofer:fig:LPTN}.
This approach is computationally efficient and often performs well after calibration to measurement data. However, LPTN models are inherently ad-hoc, require significant expert knowledge for
construction, and do not provide systematic control of approximation accuracy or spatial resolution.
These limitations can be overcome by reduced order models based on physically consistent descriptions of thermal transfer by partial differential equations (PDEs).

In this work, we consider Reduced Basis methods constructed from high-fidelity finite element approximations of the underlying PDE system. The RB framework provides systematic low-dimensional approximations of parametrized PDEs with rigorous error control and efficient online evaluation~\cite{haslhofer:hesthaven,haslhofer:quarteroni}.
The key modeling steps include:
\begin{itemize}
    \item 3D heat transfer equations in stator and rotor using co-rotating coordinate systems;
    \item a spatially homogeneous air-gap temperature representing average temperatures;
    \item appropriate interface conditions describing heat-exchange within the subdomains. 
\end{itemize}
This strategy eliminates the need for modelling rotation explicitly. The heat-transfer between subdomains is described by appropriate interface conditions which allow for parameter estimation and to include dependence on the speed-of-rotation.

All methods studied in this paper lead to parameter dependent differential equations 
\begin{align}
 \label{haslhofer:eq:LPTn}
 	C(\mu) \, \dot{u} + K(\mu) \, u = F(\mu)
\end{align}
with model parameters $\mu$ representing heat capacitances, conductances, and sources. Calibration of these parameters allows to compensate for unmodeled effects and to significantly improve the accuracy of the reduced order models.

\medskip
\noindent\textbf{Contributions and outline.}
This paper develops a reduced-order modeling framework
for transient heat transfer in rotating electric machines. We introduce a model for heat transfer across the air gap, analyze the resulting PDE model and its Galerkin discretizations (finite element and RB), and investigate parameter identification using measurement data from an induction machine. The proposed approach is assessed against FE and LPTN models in terms of accuracy, computational
efficiency, and modeling flexibility. 
The remainder of the paper is organized as follows:
Section~\ref{haslhofer:sec:2} introduces the mathematical model, 
Section~\ref{haslhofer:sec:3} describes the discretization by finite elements,
Section~\ref{haslhofer:sec:4} is concerned with the RB approximation, and 
Section~\ref{haslhofer:sec:5} presents numerical results and parameter estimation studies.


\section{Modeling Heat Transfer in Rotating Electric Machines}
\label{haslhofer:sec:2}

\begin{figure}[ht!]
\centering
\includegraphics[width=0.45\textwidth]{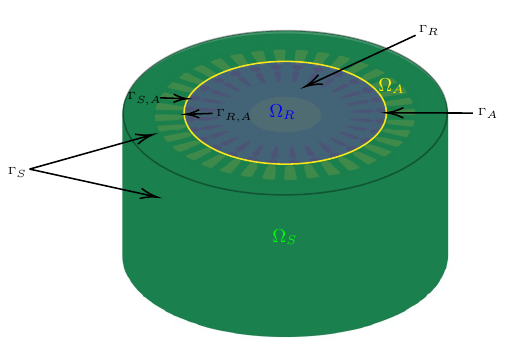}
\includegraphics[width=0.45\textwidth]{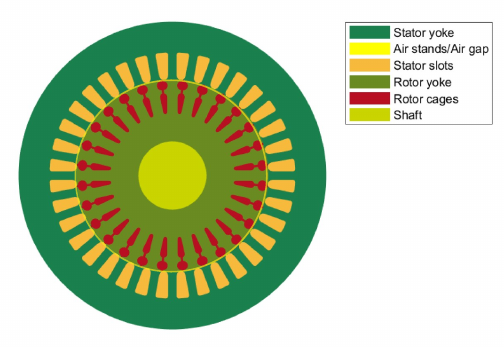}
\caption{Model geometry of the induction machine (left) and corresponding cross section (right). Colors indicate different material regions.}
\label{haslhofer:fig:fe_geo}
\end{figure}

We consider an induction machine consisting of stator, rotor, and air gap; see Figure~\ref{haslhofer:fig:fe_geo}. The heat transfer in stator and rotor is described in co-rotating coordinate systems by 
\begin{alignat}{5}
c_D \partial_t u_D
-\mathrm{div}(\kappa_D \nabla u_D)
&= f_D  \qquad &&\text{in } \Omega_D,
 \label{haslhofer:eq1}
 \\
(\kappa_D\nabla u_D)\cdot n
+
\gamma_D\, (u_D-\bar u_\infty) &= 0 
\qquad &&\text{on }\Gamma_D.
\label{haslhofer:eq2}
\end{alignat}
Here $\Omega_D$, $D \in \{S,R\}$, denotes stator and rotor domain, $\Gamma_D$ the boundary to the surrounding, $u_D$ the temperature in the stator and rotor, $c_D$ the volumetric heat capacity, $\kappa_D$ the thermal conductivity, and $f_D$ the power loss density; $\gamma_D$ is a heat transfer coefficient and $\bar u_\infty$ the ambient temperature.  
The coefficients $c_D$, $\kappa_D$, and $\gamma_D$ vary between the material regions shown in Figure~\ref{haslhofer:fig:fe_geo}. Since laminated steel sheets and windings are represented by homogenized materials, the conductivity tensor $\kappa_D$ is in general anisotropic.
Details are described in~\cite{haslhofer:msc_thesis}.

The air gap is treated as a separate subdomain $\Omega_A$. We assume that the air-gap temperature is spatially homogeneous and denote it by $\bar u_A$.  This approximation is justified by the small width of the air gap and the strong mixing induced by rotation. Heat transfer between air gap,  stator/rotor, and the surrounding is described by 
\begin{alignat}{5}
(\kappa_D\nabla u_D)\cdot n + \alpha_D \, (u_D-\bar u_A) 
&= 0 \qquad  &&\text{on } \Gamma_{D,A}, \label{haslhofer:eq3}
 \\
\bar c_A \partial_t \bar{u}_A + \bar \gamma_A (\bar{u}_A-\bar{u}_{\infty}) + \bar \alpha_S (\bar{u}_A -\bar{u}_S)+\bar \alpha_R (\bar{u}_A -\bar{u}_R)&=0 \qquad &&\text{in } \Omega_A. \label{haslhofer:eq4}
\end{alignat}
Here $\alpha_D$ is a heat transfer coefficient, and $\bar c_A$, $\bar \gamma_A$, $\bar \alpha_D$ are appropriately scaled by the air gap volume and the interface areas. Furthermore, $\bar u_D$ denotes the average of the temperature $u_D$ along the interface $\Gamma_{D,A}$ between $\Omega_D$ and $\Omega_A$. The $26$ parameters $c_D$, $\kappa_D$, $\gamma_D$ and $\alpha_D$ characterizing different materials and interfaces are collected in the parameter vector $\mu$ in the following discussion and will be estimated later on.
%

\section{Variational Formulation and Finite Element Approximation}
\label{haslhofer:sec:3}

We now formulate the thermal model of Section~\ref{haslhofer:sec:2} in a weak setting suitable for analysis and discretization. The temperature field is written as
\(u=(u_S, u_R, u_A)\). 
For the subsequent analysis, we introduce the Hilbert spaces
\[
V := H^1(\Omega_S)\times H^1(\Omega_R)\times \mathbb{R}, 
\qquad
H := L^2(\Omega_S)\times L^2(\Omega_R)\times \mathbb{R}.
\]
The weak formulation is obtained in the standard way by testing the governing equations with functions in $V$ and applying integration by parts. This leads to
\begin{problem1} \label{haslhofer:p}
Find $u \in L^\infty(0,T;V)$ with $u' \in L^2(0,T;H)$ such that
\begin{alignat}{5}
c(u'(t),v;\mu) + a(u(t),v;\mu) &= \ell(t;v;\mu) \qquad 
&& \forall v \in V,\, t\in(0,T), \\
c(u(0),v;\mu) &= c(u_0,v;\mu)
&& \forall v \in V.
\end{alignat}
\end{problem1}
Here $a(\cdot,\cdot;\mu)$ and $c(\cdot,\cdot;\mu)$ are parameter-dependent bilinear forms incorporating heat diffusion, coupling between subdomains, and thermal capacity effects, while $\ell(\cdot;\mu)$ represents heat sources and boundary contributions. The parameter vector $\mu$ collects material properties such as thermal conductivities, heat capacities, and heat transfer coefficients, and will later be used for calibration to measurement data.

\begin{assumption} \label{haslhofer:ass1}
The bilinear forms
$a(\cdot,\cdot;\mu)$ and $c(\cdot,\cdot;\mu)$ are symmetric, continuous, and
coercive on $V$ and $H$, respectively, uniformly with respect to $\mu \in \mathcal{P}$. The admissible parameter set $\mathcal{P}$ is finite dimensional and compact. 
\end{assumption}

These assumptions are easily verified in practice. They guarantee well-posedness of the continuous problem \cite{haslhofer:msc_thesis, haslhofer:evans} and will be used for error analysis below.

\bigskip 
\noindent 
\textbf{Truth approximation.}
In practice, the system is solved after spatial and temporal discretization. 
Let $V_h \subset V$ be a finite-dimensional conforming subspace obtained from a standard finite element discretization of the stator and rotor domains. The air-gap component remains scalar. All subsequent reduced-order models are constructed from this discrete reference space.
For time integration, we employ the implicit Euler scheme with uniform time step $\tau>0$. 
The resulting semi-discrete problem reads:
\begin{problem1}\label{haslhofer:discrete_p}
Find $u_h^n \in V_h$ such that for $n=0,\dots,N_t$:
\begin{alignat}{5}
c(d_\tau u_h^n, v_h;\mu) + a(u_h^n, v_h;\mu)
&= \ell(t^n; v_h;\mu) \qquad 
&& \forall v_h \in V_h, \\
c(u_h^0,v_h;\mu)
&=c(u_0,v_h;\mu)
&& \forall v_h\in V_h.
\end{alignat}
Here $d_\tau u_h^n := (u_h^n - u_h^{n-1})/\tau$ denotes the backward difference quotient.
\end{problem1}

The implementation of Problem~\ref{haslhofer:discrete_p} leads to a time-discrete version of system \eqref{haslhofer:eq:LPTn}. 
Assumption~\ref{haslhofer:ass1} allows to establish the well-posedness of the discrete problem and to show that the error between $u$ and $u_h$ can be made arbitrarily small by choosing $h$ and $\tau$ small enough; see~\cite{haslhofer:thomee}. 
This high-dimensional system serves as the reference solution for all subsequent reduced-order constructions, i.e., as a \emph{truth approximation}~\cite{haslhofer:hesthaven,haslhofer:quarteroni}.

\section{Reduced Basis Method}
\label{haslhofer:sec:4}

The Reduced Basis (RB) method constructs low-dimensional approximations of the truth approximation by restricting the discrete problem
to a suitably chosen subspace
\[
V_N \subset V_h, \qquad N := \dim(V_N) \ll \dim(V_h).
\]
A low dimensional approximation of Problem~\ref{haslhofer:discrete_p} is then obtained by restriction to $V_N$.
\begin{problem1}\label{haslhofer:reduced_p}
Find $u_N^{n} \in V_N$ for $n =0,\ldots,N_t$, such that 
\begin{align} 
	c(d_\tau u_N^{n},v_N;\mu) +a(u_N^{n},v_N;\mu)&=\ell(t^{n};v_N;\mu) &\qquad \forall v_N \in V_N,\label{haslhofer:reduced_eq1} \\
	c(u_N^{0},v_N;\mu)&=c(u_0,v_N;\mu) &\qquad \forall v_N \in V_N. \label{haslhofer:reduced_eq2} 
\end{align}
\end{problem1}
The two Problems~\ref{haslhofer:discrete_p} and \ref{haslhofer:reduced_p} lead to time-discrete versions of system \eqref{haslhofer:eq:LPTn} with dimension $n$ and $N$, respectively. The main goal of the RB approximation is therefore to construct $V_N$ such that
$N$ remains small while maintaining high fidelity with respect to the truth approximation.
This can be motivates the following considerations.

\medskip
\noindent
\textbf{Error structure of the RB approximation.}
The following result shows that the RB error can be fully controlled by approximation properties of the reduced space with respect to the FE solution
manifold; see the discussion below.

\begin{theorem1}
\label{haslhofer:thm:main}
Let Assumption~\ref{haslhofer:ass1} hold and $\vertiii{u_h}^2 = \sum_n \tau \|u_h^n\|^2_V$. Then 
\[
\vertiii{u_h(\mu)-u_N(\mu)}^2
\le
C 
\vertiii{u_h(\mu)-R_N u_h(\mu)}^2
\]
where $R_N:V\to V_N$ denotes the Ritz projection based on the bilinear form $a(\cdot,\cdot;\mu)$.
Moreover, the constant $C$ is independent of $V_h$, $V_N$, $\tau$, and $\mu \in \mathcal{P}$.
\end{theorem1}

The proof follows standard arguments for Galerkin approximations of parabolic problems; see~\cite{haslhofer:thomee}. For detailed derivation of similar estimates, we refer to \cite{haslhofer:msc_thesis}. 

\medskip
\noindent
\textbf{Interpretation.}
Theorem~\ref{haslhofer:thm:main} shows that the RB error is fully determined by how well the finite element snapshots $u_h^n(\mu)$ can be approximated in the reduced space $V_N$. In particular, uniform accuracy follows if
\[
\|u_h^n(\mu) - R_N u_h^n(\mu)\|_V \le \varepsilon
\qquad \forall \mu \in \mathcal{P},\; n=1,\dots,N_t.
\]
This motivates the definition of the parameter-dependent solution manifold
\[
\mathcal{M}_h :=
\{ u_h^n(\mu) : \mu \in \mathcal{P},\; 0 \le t^n \le T \}
\subset V_h,
\]
and the construction of $V_N$ can be interpreted as seeking a low-dimensional space that approximates $\mathcal{M}_h$ efficiently.
Standard results in reduced basis theory \cite{haslhofer:quarteroni} show that low-dimensional approximation is possible when the solution depends sufficiently smoothly on the parameters and the underlying problem data. In this case, $\mathcal{M}_h$ can often be approximated with rapidly decaying error using small-dimensional spaces $V_N$.

\bigskip
\noindent 
\textbf{Construction of the reduced basis.}
The setup of the RB space follows standard arguments~\cite{haslhofer:quarteroni}. For completeness of the presentation, we provide some details about the practical implementation we use for our numerical tests later on. 
\begin{enumerate}
\item \textbf{Snapshot generation:}
    The full-order problem is solved for selected parameter values $\mu_j \in \mathcal{P}$, $j=1,\dots,n_j$. The resulting trajectories $u_h^n(\mu_j)$ are collected as columns of the snapshot matrix $S$. Parameter samples are generated by random perturbations around reference configurations as described in~\cite{haslhofer:msc_thesis}. Comparable results are obtained with other sampling strategies, e.g., Latin hypercube sampling.

\item \textbf{Basis extraction:}
    A reduced basis is obtained using Proper Orthogonal Decomposition (POD). To alleviate the computational cost of the SVD of the large snapshot matrix, we employ a randomized SVD strategy as proposed in~\cite{haslhofer:tropp}. The reduced basis consists of the first $N$ left singular vectors, where $N$ is chosen such that
    \[
    \frac{\sum_{i=1}^{N} \sigma_i^2}{\sum_{i=1}^{r} \sigma_i^2}
    \ge 1 - \varepsilon_{\mathrm{POD}}^2,
    \qquad \sigma_1 \ge \sigma_2 \ge \dots \ge \sigma_r,
    \]
    with $r$ denoting the rank of the snapshot matrix. 
    In accordance with Theorem~\ref{haslhofer:thm:main}, the POD is performed with respect to the energy-norm $\|\cdot\|_V$; compare with~\cite{haslhofer:quarteroni}. 

\item \textbf{Reduced basis space:}
    The selected POD modes define the reduced space $V_N := \mathrm{span}\{\phi_1,\dots,\phi_N\} \subset V_h$. This step yields the discrete trial space used in the Galerkin projection. This space is represented algebraically by projection matrix, which allows to assemble the main building blocks of the RB method in the \emph{offline stage}.

\item \textbf{Affine decomposition:}
    The particular modeling approach for the air-gap coupling provides affine parameter dependence~\cite{haslhofer:quarteroni}; see eq.~\eqref{haslhofer:eq1}--\eqref{haslhofer:eq4} and also \cite{haslhofer:msc_thesis} for details. The system matrices of the RB approximation appearing in \eqref{haslhofer:eq:LPTn} thus have the form 
    \begin{align*}
    K(\mu) = \sum\nolimits_q \theta^q_a(\mu) K^q, \quad 
    C(\mu) = \sum\nolimits_q \theta^q_c(\mu) C^q, \quad \text{and} \quad 
    F(\mu) = \sum\nolimits_q \theta^q_f(\mu) F^q.
    \end{align*}
    Since the components $K^q$, $C^q$, $F^q$ can be pre-computed in the offline phase, the assembly cost on the \emph{online phase} becomes negligible. 

\item \textbf{Online evaluation:}
   The RB approximation leads to systems of the form \eqref{haslhofer:fig:LPTN} of small size $N$. Since the systems are time-invariant, we use diagonalization of the matrices to allow for efficient time-integration. Assuming $N \ll N_t$, the online-complexity of a forward simulation can thus be reduced to $O(N \cdot N_t)$ operations, which is optimal.
\end{enumerate}
\medskip
\noindent
The resulting reduced basis model provides a fully computable low-dimensional approximation of the thermal system, combining data-driven basis construction with structure-preserving projection of the underlying PDE model. The affine decomposition and the smooth dependence of the solution on $\mu$ provide reducibility of the problem and are therefore key ingredients of the RB method, see~\cite{haslhofer:quarteroni}. In the following section, we investigate the numerical performance of the proposed approach. In particular, we compare the RB method with the full finite element model and classical Lumped Parameter Thermal Networks in terms of accuracy, computational cost, and robustness with respect to parameter variations and measurement-based calibration.

\section{Numerical validation}
\label{haslhofer:sec:5}

\begin{figure}[ht!]
	\includegraphics[ height=0.43\textwidth]{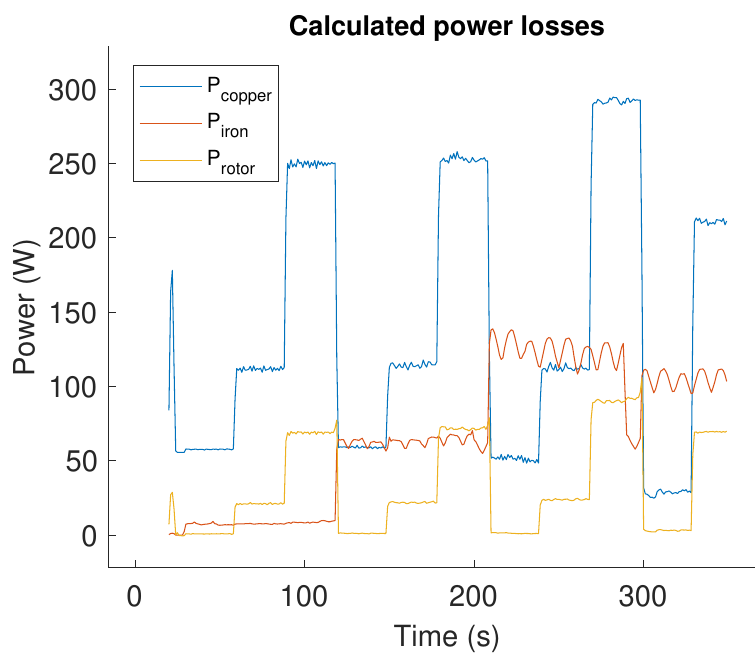}
   \hfill
     \includegraphics[height=0.44\textwidth]{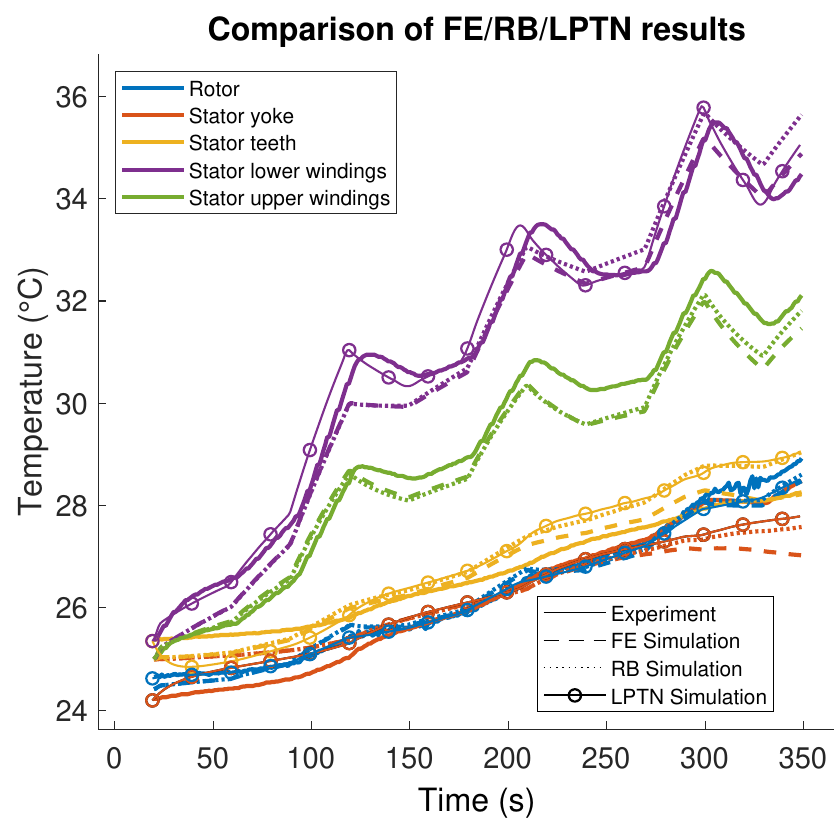}
\caption{Calculated power losses for conducted experiment (Left)~\cite{haslhofer:heidarikani} and comparison of  FE/RB/LPTN results in sensor points with measurements (Right).}%
\label{haslhofer:fig:results}%
\end{figure}

We now compare the accuracy and performance of finite element (FE), reduced basis (RB), and lumped parameter thermal network (LPTN) for reproducing measurements from a real induction machine experiment conducted at TU Graz~\cite{haslhofer:heidarikani}. 
In the experiment, three-phase AC is supplied such that the rotational speed is increased iteratively from 0 to 2500 rpm, while the motor load is increased during each speed interval. Temperature sensors detect the temperatures in the stator and rotor and power losses are calculated using the measured input power, frequency and flux linkage, together with analytical formulas~\cite{haslhofer:heidarikani}, see Figure~\ref{haslhofer:fig:results}. 
To obtain good agreements with measurements, parameter estimation is performed for all computational models using a nonlinear least-squares approach. All computations are carried out
in \textsc{Matlab} and on a desktop computer with 11th Gen Intel Core i7-1165G7 CPU, clock rate 2.701 GHz.
The results of our validation studies are reported in Figure~\ref{haslhofer:fig:results}. 

\medskip
\noindent
\textbf{Implementation details.}
To ensure a fair comparison between all methods, several common and method-specific implementation choices are made. For the FE model, geometric symmetries are exploited such that only one stator slot and one rotor bar are resolved; this reduction is justified by the structure of the air-gap coupling. In addition, a tensor-product discretization is used, combining a triangular finite element mesh in the cross-section with global polynomial basis functions in the axial direction. Since temperature variations in axial direction are small, only three axial modes are retained, resulting in a highly efficient FE approximation.
All methods employ the same implicit Euler time-stepping scheme to ensure consistency of the temporal discretization with a uniform time step $\tau=1\,s$. For the RB and LPTN models, system matrices are diagonalized, leading to optimal online complexity of order $O(N \times N_t)$. For the FE model, a pre-factorization of the system matrices is used, yielding an efficient online evaluation with complexity close to $O(n \log(n)\, N_t)$, where $n$ denotes the FE dimension.

\noindent 
\textbf{Results.} 
All considered methods achieve comparable accuracy after calibration. 
%
%
The system dimension and thus computational efficiency, however,  differs significantly between the high-fidelity method (FE: $n=1032$) and the reduced basis approximations (RB:  $N=5$;  LPTN: $N=5$). The LPTN and RB Method not only decrease computation times but also significantly reduce the memory requirements; see Table~\ref{haslhofer:tab:results} for details. 

\begin{table}[ht!]
\centering
\renewcommand{\arraystretch}{1.3}
\setlength{\tabcolsep}{10pt}
\begin{tabular}{l||r|r|r}
     & FE Model & RB Model & LPTN \\
  \hline
  \hline
  Computation time & 0.013 s & 0.004\,s & 0.004\,s \\
  \hline
  Online memory    & 84\,KB    & 200\,B    & 200\,B \\
\end{tabular}
\caption{Comparison of FE, RB, and LPTN models in terms of computational performance for a single simulation, including both execution time and memory requirements. The results highlight the efficiency gains of reduced-order approaches compared to the full finite element model.}
\label{haslhofer:tab:results}
\end{table}
Both RB and LPTN models achieve a substantial reduction in computational ressources, while providing similar accuracy to FE simulations, which seems particularly relevant for real-time embedded applications.
The RB approach further offers a systematic and fully automated construction of reduced-order spaces together with rigorous control of the approximation quality. This is in contrast to LPTN models, which typically rely on manual model design. Moreover, RB models provide information about spatial temperature distributions, while LPTNs only yield compartment averages. 

\begin{acknowledgement}
This research was funded by the Austrian Science Fund (FWF) 10.55776/F90.
\end{acknowledgement}

\end{document}